\input amstex
\documentstyle{amsppt}

\magnification =\magstep1
\refstyle {C}

\def\dd       #1#2#3{#1_{#2#3}}
\def\ddd      #1#2#3#4{#1_{#2#3#4}}
\def\dddd     #1#2#3#4#5{#1_{#2#3#4#5}}

\def\ud       #1#2#3{{#1}^{#2}_{\phantom{#2}{#3}}}
\def\du       #1#2#3{{#1}_{#2}^{\phantom{#2}{#3}}}

\def\ddu      #1#2#3#4{{#1}_{#2#3}^{\phantom{#2#3}{#4}}}
\def\udd      #1#2#3#4{#1^#2_{\phantom{#2}#3#4}}
\def\duu      #1#2#3#4{#1_#2^{\phantom{#2}#3#4}}

\def\uddd     #1#2#3#4#5{#1^#2_{\phantom{#2}#3#4#5}}
\def\dudd     #1#2#3#4#5{#1_{#2\phantom{#3}#4#5}^{\phantom{#2}#3}}

\def\dddu     #1#2#3#4#5{#1_{#2#3#4}^{\phantom{#2#3#4}#5}}

\def\pdu      #1#2#3#4{{\partial {#1}^{#2}\over\partial {#3}^{#4}}}

\def\pddu     #1#2#3#4#5#6{{\partial^2{#1}^{#2}
              \over\partial {#3}^{#4}\partial {#5}^{#6}}}

\def\td       #1#2#3{{d {#1}\over d {#2}_{#3}}}

\def\cd       #1#2{\nabla_{#2}{#1}}
\def\cdd      #1#2#3{\nabla_#2\nabla_#3#1}

\def\cdb      #1#2{\overline{\nabla}_{#2}{#1}}
\def\cddb     #1#2#3{\overline{\nabla}_{#2}\overline{\nabla}_{#3}{#1}}

\def\con      #1#2#3{\Gamma^{#1}_{#2#3}}
\def\conbar   #1#2#3{\overline{\Gamma}^{#1}_{#2#3}}

\def\contilde #1#2#3{\widetilde{\Gamma}^{#1}_{#2#3}}

\def\t        #1{\widetilde {#1}}

\def\b        #1{\overline{#1}}

\def\formula  #1#2{$$#2\leqno#1$$}
\def\dt       {{d\over dt}}


\topmatter
\title        A canonical way to deform a Lagrangian submanifold
\endtitle
\author       Knut Smoczyk
\endauthor
\address      Mathematics Department, Harvard University, Cambridge, MA
              02138, USA
\endaddress
\email        ksmoczyk\@abel.math.harvard.edu
\endemail
\date         May 22, 1996
\enddate
\thanks       Supported by Alexander von Humboldt Foundation
\endthanks
\subjclass    53C42, 53C25, 53C55
\endsubjclass
\abstract     We derive some important geometric identities for Lagrangian
submanifolds immersed in a K\"ahler manifold and prove that there exists
a canonical way to deform a Lagrangian submanifold by a parabolic flow
through a family of Lagrangian submanifolds if the ambient space is a
Ricci-flat Calabi-Yau manifold.
\endabstract
\endtopmatter
\document

\baselineskip =12pt 

\head         0. Introduction
\endhead

Let $~L^n~$ be a smooth submanifold, immersed into a Ricci-flat Calabi-Yau
manifold $~M^{2n}~$ with complex structure $~J~$ and metric $~\b g~$. Let
$~\b\omega=\b g(J\cdot,\cdot)~$ be the K\"ahlerform on $~M~$ and let
$~g,~\omega~$ denote the pullbacks on $~L~$. If $~\omega=0~$, then $~L~$
is called Lagrangian. It is an interesting question,
if there exists a canonical way to deform an initial Lagrangian submanifold
$~L_0~$
through a family of submanifolds $~L_t~$ such that each $~L_t~$ is also
Lagrangian, i.e. $~\dt \omega=0~$ and such that
the corresponding flow is parabolic.
Given a 1-form on a Lagrangian submanifold, one can first use the metric
$~g~$ to identify
this 1-form with a vectorfield and then apply the
complex structure which by the Lagrangian condition maps this tangent
vector field
to a normal vector field. Assuming that we deform the Lagrangian submanifold
in this direction one obtains the necessary condition that this 1-form
has to be closed in order to maintain the Lagrangian structure. This has
been proven in \cite 6 but will also become clear in this paper.
However, it is not obvious if this is also sufficient.
So if one wants to find such a canonical deformation
the first thing one has to make sure is that there exists
a canonical, closed 1-form on a Lagrangian submanifold.
Surprisingly this is true if the ambient space is Ricci-flat. The
corresponding 1-form can be defined in terms of the second fundamental
form and we denote this 1-form by ``mean curvature form" (see definition
below) since the resulting deformation vector field is given by the mean
curvature vector which can be defined for arbitrary smooth submanifolds.
The organization of the paper is as follows:

In paragraph 1 we state important geometric equalities
some of which are analogous to the Gauss-Weingarten, Codazzi equations, etc.
and define the mean curvature form. Then we prove that the
corresponding flow exists at least on a short time interval and has a
unique solution on the moduli space of smooth submanifolds. Finally we derive
an evolution inequality
for the squared norm of $~\omega~$ which basically follows from the
Ricci-flatness of the ambient space and with a parabolic maximum principle
we can prove that $~\omega~$ has to vanish identically if
it vanishes at $~t=0~$. Thus, if one wants to study the geometry of Lagrangian
submanifolds immersed into a Ricci-flat Calabi-Yau manifold this flow
provides a tool and there should be no doubt that this is a powerful
tool since the corresponding mean curvature flow for hypersurfaces in
Riemannian manifolds has already given many results (e.g. \cite 5).

In paragraph 2 we restate some equations derived in paragraph 1 which in
the context of a Lagrangian submanifold become very beautiful and
eventually derive the evolution equation for the second fundamental form and
the mean curvature form. We also state the evolution equations in the
case where $~L~$ is being deformed by an arbitrary 1-form (formally).
As a direct consequence of these evolution equations we finally
prove a useful theorem and the paper ends with some questions
and remarks.

\head           1. Notations and preliminary results
\endhead

Let us define the following map:
$$~N:T_pL^n\to (T_pL^n)^\perp~,~N(v):=(J(v))^\top~,$$
where $~(T_pL^n)^\perp~$ is the normal space of
$~L~$ at $~p~$ and $~^\top~$ denotes
the projection onto $~T_pL^n~$. In the forthcoming let us assume that $~L~$ is
a compact submanifold such that $~N~$ is always an isomorphism between $~T_pL~$
and $~(T_pL)^\perp~$. Then we define the following tensor
on $~L~$:
$$h(u,v,w):=-\b g(N(u),\cdb wv)=\b g(\cdb {N(u)}v,w)$$
where $~\b\nabla~$ denotes the covariant derivative on $~M~$. For a fixed
vector $~u~$ this is the second fundamental form with respect to the normal
vector $~N(u)~$.
Assume that $~F:L^n\to M^{2n}~$ is an immersion, that $~x^i,~y^\alpha
,~i=1\dots n,~\alpha=1\dots 2n~$
are coordinates for $~L^n,~M^{2n}~$ respectively
and set $~e_i:={\partial F^\alpha\over\partial x^i}{\partial\over\partial
y^\alpha}~$, where double indices are always summed from $~1~$ to $~n~$ or
to $~2n~$ respectively. Furtheron we will always write $~\langle u,v\rangle~$
instead of $~\b g(u,v)~$. Then
$$\eqalign{g_{ij}&=\b g_{\alpha\beta}\pdu F{\alpha}xi\pdu F{\beta}xj\cr
N(e_i)&=J(e_i)-\du {\omega}ile_l\cr
\b\omega_{\alpha\beta}&=\b g_{\beta\gamma}\ud J{\gamma}{\alpha}\cr
\b\omega_{\alpha\beta}&=-\b\omega_{\beta\alpha}\cr
\b\nabla J&=0\cr}$$
Since $~h(u,\cdot,\cdot)~$ is the second fundamental form with respect to
$~N(u)~$, we clearly have
$$h(e_k,e_i,e_j)=h_{kij}=h_{kji}$$
Now we ask: Under what condition
is $~h_{ijk}~$ also symmetric in the first two indices?
The answer is

\proclaim {Proposition 1.1}
$~h_{kij}=h_{ikj}+\cd {\omega_{ik}}{j}~$
\endproclaim
\demo{Proof}
$$\eqalign{h_{kij}=h_{kji}
&=-\langle N(e_k),\cdb {e_i}{e_j}\rangle=-\langle J(e_k),
\cdb {e_i}{e_j}\rangle+\du {\omega}kl\langle e_l,\cdb {e_i}{e_j}\rangle\cr
&=\langle e_k,\cdb {J(e_i)}{e_j}\rangle+\du {\omega}kl\langle e_l,
\cdb {e_i}{e_j}\rangle\cr
&=\omega_{ik,j}-\langle N(e_i),\cdb {e_k}{e_j}\rangle-\du {\omega}il\langle
e_l,\cdb {e_k}{e_j}\rangle+\du {\omega}kl\langle e_l,\cdb {e_i}{e_j}\rangle\cr
&=h_{ijk}+\omega_{ik,j}-\du {\omega}il\langle e_l,(\cdb {e_k}{e_j})^\top
\rangle-\ud {\omega}lk\langle e_l,(\cdb {e_i}{e_j})^\top\rangle\cr
&=h_{ikj}+\cd {\omega_{ik}}j\cr}$$
where we have used that
$$\eqalign{\langle J(u),J(v)\rangle&=\langle u,v\rangle\cr
\b\nabla J&=0\cr
J^2&=-id\cr}$$
\enddemo

We will need the following identity
\proclaim{Proposition 1.2}
$$\cdd {\omega_{li}}kj-\cdd {\omega_{ki}}lj=\cdd {\omega_{lk}}ji+\uddd Rsilj
\omega_{ks}+\uddd Rsijk\omega_{ls}-\uddd Rsjkl\omega_{si}$$
\endproclaim
\demo{Proof}
This follows from the rule for interchanging derivatives and the fact that
$~\omega~$ is closed:
$$\eqalign{\cdd {\omega_{li}}kj&=\cdd {\omega_{li}}jk
+\uddd Rsijk\omega_{ls}+\uddd Rsljk\omega_{si}\cr
&=-\cdd {\omega_{ik}}jl-\cdd {\omega_{kl}}ji+\uddd Rsijk\omega_{ls}
+\uddd Rsljk\omega_{si}\cr
&=\cdd {\omega_{ki}}lj+\uddd Rsilj\omega_{ks}-\uddd Rsjkl\omega_{si}
-\cdd {\omega_{kl}}ji+\uddd Rsijk\omega_{ls}\cr}$$
where we have used the Bianchi identity in the last step.
\enddemo

\remark {Remark}
If $~v_{\alpha},~\alpha =1\dots 2n~$ is a basis for $~T_pM^{2n}~$ and
$~a^{\alpha\beta}~$ is the inverse of
$~a_{\alpha\beta}:=\langle v_{\alpha},v_{\beta}\rangle~$, then any vector
$~w\in T_pM^{2n}~$ can be written in the form $~w=a^{\alpha\beta}\langle
w, v_{\alpha}\rangle v_{\beta}~$.
\endremark

Now define the following tensor
$$\eta_{ij}:=\langle N(e_i),N(e_j)\rangle=g_{ij}+\du {\omega}il\dd {\omega}lj$$
This is a positive definite, symmetric tensor as long as $~N~$ is an
isomorphism. Let us also denote the inverse of $~\eta_{ij}~$ by $~\eta^{ij}~$.
In the forthcoming we will always assume that $~e_i,~i=1\dots n~$ are
not only defined on $~L~$ but that we have chosen an extension in a
tubular neighborhood of $~L~$. All the stated equations like
$~\cdb {e_j}{e_i}=\dots~$, etc. will mean that this is true if we
restrict these terms to the submanifold $~L~$. In particular it will
always be true that these results are independent of the chosen extension
for $~e_i,~i=1\dots n~$.

\proclaim{Proposition 1.3}
\formula{(a)}{\cdb {e_j}{e_k}=\con nkje_n-\eta^{mn}h_{mkj}N(e_n)}
$$\leqalignno{\cd {h_{ikj}}l-\cd {h_{ilj}}k&=\dddd {\b R}{\underline i}jkl
+\eta^{mn}\du{\omega }ns (h_{mlj}h_{ski}-h_{mkj}h_{sli})&(b)\cr
&+\eta^{mn}\du{\omega }is (h_{mkj}h_{nls}-h_{mlj}h_{nks})& \cr}$$
\formula {(c)}{\dddd Rijkl=\dddd {\b R}ijkl+\eta^{mn}(h_{mik}h_{njl}
-h_{mil}h_{njk})}
\endproclaim
where an underlined index means that one has to take the image of this vector
under $~N~$, e.g.
$~\dddd {\b R}{\underline i}jkl=\langle \b R(e_k,e_l)e_j,N(e_i)\rangle~$

\demo{Proof}

(a): The above remark with $~v_i=e_i,~v_{i+n}=N(e_i),~i=1\dots n~$ gives
$$\eqalign{\cdb {e_j}{e_k}&=g^{mn}\langle\cdb {e_j}{e_k},e_m\rangle e_n+
\eta^{mn}\langle \cdb {e_j}{e_k},N(e_m)\rangle N(e_n)\cr
&=g^{mn}\langle\cd {e_j}{e_k},e_m\rangle e_n-\eta^{mn}h_{mkj}N(e_n)\cr
&=\con nkje_n-\eta^{mn}h_{mkj}N(e_n)\cr}$$

(b): Let us choose normal coordinates for $~L~$ at a fixed point $~p~$, then
$$~\con kij,~\cd {e_j}{e_i},~[e_i,e_j]=0,~\cd {h_{ikj}}l=h_{ikj,l}~$$
and
$$\eqalign{\cd {h_{ikj}}l-\cd {h_{ilj}}k
&=\langle\cdb {N(e_i)}{e_k},\cdb {e_j}{e_l}
\rangle-\langle\cdb {N(e_i)}{e_l},\cdb {e_j}{e_k}\rangle\cr
&+\langle N(e_i),\cddb {e_j}{e_k}{e_l}-\cddb {e_j}{e_l}{e_k}\rangle\cr
&=\langle\cdb {N(e_i)}{e_k},\cdb {e_j}{e_l}
\rangle-\langle\cdb {N(e_i)}{e_l},\cdb {e_j}{e_k}\rangle\cr
&+\langle N(e_i),\b R(e_k,e_l)e_j\rangle\cr
}$$
With (a) we get
$$\cd {h_{ikj}}l-\cd {h_{ilj}}k=\dddd {\b R}{\underline i}jkl-\eta^{mn}\bigl(
h_{mlj}\langle \cdb {N(e_i)}{e_k},N(e_n)\rangle
-h_{mkj}\langle \cdb {N(e_i)}{e_l},N(e_n)\rangle\bigr)$$
Now we compute
$$\eqalign{\langle \cdb {N(e_i)}{e_k},N(e_n)\rangle&=\langle \cdb {\bigl(
J(e_i)-\du {\omega}is e_s\bigr)}{e_k},N(e_n)\rangle\cr
&=\langle \cdb {J(e_i)}{e_k},N(e_n)\rangle-\du {\omega}is\langle
\cdb {e_s}{e_k},N(e_n)\rangle\cr}$$
since $~\langle e_s,N(e_n)\rangle=0~$. Then
$$\eqalign{\langle \cdb {N(e_i)}{e_k},N(e_n)\rangle
&=-\langle\cdb {e_i}{e_k},J(N(e_n))\rangle+\du {\omega}ish_{nks}\cr
&=\langle\cdb {e_i}{e_k},e_n+\du{\omega}ns\du {\omega}spe_p+
\du {\omega}nsN(e_s)\rangle+\du {\omega}ish_{nks}\cr
&=-\du {\omega}nsh_{ski}+\du {\omega}ish_{nks}\cr}$$
Inserting this equality in the above identity gives the desired result.

(c): This is only a reformulation of the Gauss curvature equations.
\enddemo

\proclaim {Lemma 1.4}
$$\eqalign{\cd {h_{kij}}l-\cd {h_{lij}}k&=\dddd {\b R}{\underline i}jkl+
\cdd {\omega_{lk}}ji+\du {\omega}is\dddd {\b R}sjkl\cr
&+\du {\omega}ks\dddd Rsilj+\du {\omega}ls\dddd Rsijk
+\eta^{mn}\du {\omega}ns(h_{mlj}h_{ski}-h_{mkj}h_{sli})\cr}$$
\endproclaim

\demo{Proof}
Using Proposition 1.1 we obtain
$$\eqalign{\cd {h_{kij}}l-\cd {h_{lij}}k&=\cd {(h_{ikj}+\cd {\omega_{ik}}j)}l
-\cd {(h_{ilj}+\cd {\omega_{il}}j)}k\cr
&=\cd {h_{ikj}}l-\cd {h_{ilj}}k+\cdd {\omega_{li}}kj-\cdd {\omega_{ki}}lj\cr}$$
and with Proposition 1.2 and 1.3 (b)
$$\eqalign{\cd {h_{kij}}l-\cd {h_{lij}}k&=
\dddd {\b R}{\underline i}jkl
+\eta^{mn}\du{\omega }ns (h_{mlj}h_{ski}-h_{mkj}h_{sli})\cr
&+\eta^{mn}\du{\omega }is (h_{mkj}h_{nls}-h_{mlj}h_{nks})
+\cdd {\omega_{lk}}ji\cr
&+\uddd Rsilj\omega_{ks}+\uddd Rsijk\omega_{ls}-\uddd Rsjkl\omega_{si}\cr}$$
and by the Gauss curvature equations we have
$$\eta^{mn}\du{\omega }is (h_{mkj}h_{nls}-h_{mlj}h_{nks})
-\uddd Rsjkl\omega_{si}=\du {\omega}is\dddd {\b R}sjkl$$
which proves the lemma.
\enddemo

\definition {Definition 1.5}
$$~H:=H_idx^i:=g^{kl}h_{ikl}dx^i$$
is called the mean curvature form on $~L~$.
\enddefinition

\remark {Remark}
Using $~\eta^{ij}~$ to raise indices we can identify this 1-form with the
vector field $~\eta^{mn}H_me_n~$ and $~N~$ maps this vector field to
the outward pointing mean curvature vector field on $~L~$
(the notion ``outward pointing" does
not really make sense in arbitrary codimensions but we will denote this
vector by outward pointing since the construction is analogous to
the construction of the outward pointing
mean curvature vector in codimension 1).
Note that if $~L~$ is
Lagrangian then $~\eta_{ij}=g_{ij}~$.
\endremark

We are now going to prove that the following evolution exists at least on some
short time interval and is unique if $~L~$ is compact.
Therefore assume that $~L~$ is a compact manifold
smoothly immersed into a Ricci-flat Calabi-Yau manifold by
an immersion $~F~$ and assume that $~\eta_{ij}~$ is positive definite
everywhere on $~L~$. We want to find a smooth family of immersions
$~F_t~$ such that
\formula {(*)}{\dt F_t=-\eta^{mn}H_mN(e_n)~,~F_0=F}

\proclaim {Proposition 1.6}
$(*)$ has a unique solution on a short time interval.
\endproclaim

\demo {Proof}
The proof is formally the same as for the mean curvature flow in codimension 1.
We will use Hamilton's existence theorem for evolution equations with an
integrability condition \cite 2. Let $~\pi_N~$ be the projection on the
normal bundle
and $~\pi_L~$ be the projection on the tangent bundle of $~L~$. Here we have
$~\dt F=E(F)=-\eta^{mn}H_mN(e_n)~$. On the symbol level we get
$$\eqalign{\t g_{ij}&=\b g_{\alpha\beta}\bigl(\pdu F\alpha xi\t F^\beta \xi_j
+\pdu F\beta xj\t F^\alpha \xi_i\bigr)\cr
\contilde lij&={1\over 2}g^{kl}
(\t g_{jk}\xi_i+\t g_{ik}\xi_j-\t g_{ij}\xi_k)\cr}$$
which combine to
$$\contilde lij=g^{kl}\b g_{\beta\gamma}\pdu F\beta xk \t F^\gamma\xi_i\xi_j$$
If we multiply 1.3 (a) with $~g^{kj}~$ we see that the symbol is given by
$$\eqalign{\sigma DE(F)\t F^\alpha(\xi)&=g^{ij}(\xi_i\xi_j
\t F^\alpha- g^{kl}\b g_{\beta\gamma}\pdu F\beta xk
\t F^\gamma\xi_i\xi_j\pdu F\alpha xl)\cr
&=\vert\xi\vert^2_g(\t F^\alpha-\pi_L(\t F)^\alpha)
=\vert\xi\vert^2_g\pi_N(\t F)^\alpha\cr}$$
Therefore we are done if we use
$~\pi_L(\t F)=0~$ as our integrability condition.
\enddemo

We will need the evolution equations for $~g_{ij}~$ and
$~\omega_{ij}~$.
\proclaim {Lemma 1.7}
\formula{(a)}{\dt g_{ij}=-2\eta^{kl}H_kh_{lij}}
\formula{(b)}{\dt d\mu=-\eta^{mn}H_mH_nd\mu}
\formula{(c)}{\dt \omega=dH}
where $~d\mu~$ is the volume form on $~L_t~$.
\endproclaim

\demo{Proof}
{}From now on let us always assume that we have chosen normal coordinates
at a fixed point $~p\in L~$ with
respect to the metric $~g_{ij}~$ such that in addition $~e_i,~i=1,\dots,n~$
are eigenvectors for $~\eta_{ij}~$, i.e. $~\eta_{ij}~$ becomes diagonal in
this coordinate frame. We will also assume that we have chosen
normal coordinates at $~p~$ for the ambient manifold $~M~$. Under this
assumptions it is easy to see that
$$\eqalign{\dt g_{ij}&=
\dt\bigl( \b g_{\alpha\beta}\pdu F{\alpha}xi\pdu F{\beta}xj\bigr)
=\b g_{\alpha\beta}{\partial \over\partial x^i}(\td {F^\alpha}t{})\pdu F{
\beta}xj+
\b g_{\alpha\beta}{\partial \over\partial x^j}(\td {F^\beta}t{})\pdu F{
\alpha}xi\cr
&=-\eta^{mn}H_m\langle \cdb {N(e_n)}{e_i},e_j\rangle
-\eta^{mn}H_m\langle \cdb {N(e_n)}{e_j},e_i\rangle\cr
&=-2\eta^{mn}H_mh_{nij}\cr}$$
which gives (a) and (b).
For $~\omega~$ we can proceed similarly
$$\eqalign{\dt \omega_{ij}&=
\dt\bigl( \b \omega_{\alpha\beta}\pdu F{\alpha}xi
\pdu F{\beta}xj\bigr)\cr
&=\b \omega_{\alpha\beta,\gamma}\td {F^\gamma}t{}
\pdu F{\alpha}xi\pdu F{\beta}xj+\b \omega
_{\alpha\beta}{\partial \over\partial x^i}
(\td {F^\alpha}t{})\pdu F{\beta}xj+\b \omega
_{\alpha\beta}{\partial \over\partial x^j}
(\td {F^\beta}t{})\pdu F{\alpha}xi\cr}$$
Then we use the closedness of $~\b \omega~$ to obtain
$$\eqalign{\dt \omega_{ij}&=-
\b \omega_{\beta\gamma,\alpha}
\td {F^\gamma}t{}\pdu F{\alpha}xi\pdu F{\beta}xj
-\b \omega_{\gamma\alpha,\beta}\td {F^\gamma}t{}
\pdu F{\alpha}xi\pdu F{\beta}xj\cr
&+\b \omega
_{\gamma\beta}{\partial \over\partial x^i}
(\td {F^\gamma}t{})\pdu F{\beta}xj+\b \omega
_{\alpha\gamma}{\partial \over\partial x^j}
(\td {F^\gamma}t{})\pdu F{\alpha}xi\cr
&={\partial\over\partial x^i}
\bigl(\b\omega_{\gamma\beta}\td {F^\gamma}t{}\pdu F{\beta}xj\bigr)
-\b\omega_{\gamma\beta}\td {F^\gamma}t{}
{\partial^2 F^{\beta}\over\partial x^j\partial x^i}\cr
&-{\partial\over\partial x^j}
\bigl(\b\omega_{\gamma\alpha}\td {F^\gamma}t{}\pdu F{\alpha}xi\bigr)
+\b\omega_{\gamma\alpha}\td {F^\gamma}t{}
{\partial^2 F^{\alpha}\over\partial x^i\partial x^j}\cr
&=
{\partial\over\partial x^i}
\bigl(\b\omega_{\gamma\beta}\td {F^\gamma}t{}\pdu F{\beta}xj\bigr)
-{\partial\over\partial x^j}
\bigl(\b\omega_{\gamma\alpha}\td {F^\gamma}t{}\pdu F{\alpha}xi\bigr)
}$$
by the antisymmetry of $~\b\omega~$. Finally
$$\eqalign{\b\omega(\td Ft{},\pdu F{}xi)&=-\eta^{kl}H_k\b\omega(N(e_l),e_i)
=-\eta^{kl}H_k\b g(J(N(e_l)),e_i)\cr
&=\eta^{kl}H_k\b g(e_l+\du {\omega}lnJ(e_n),e_i)=
\eta^{kl}H_k(g_{li}+\du {\omega}ln\dd {\omega}ni)\cr
&=\eta^{kl}H_k\eta_{li}=H_i\cr}$$
and therefore
$$\dt \omega_{ij}=H_{j,i}-H_{i,j}$$
which proves (c).
\enddemo

\proclaim {Proposition 1.8}
For each compact time interval $~[0,T]~$ on which
a smooth solution of $(*)$ exists and
the flow is well defined, i.e.
where $~N~$ is an isomorphism
we can find a positive constant such that
$$\dt\vert{\omega}\vert^2\le \Delta\vert {\omega}\vert^2+c\vert{\omega}
\vert^2$$
\endproclaim

\demo{Proof}
First we will need the evolution equation for $~\vert {\omega}\vert^2
:=g^{ik}g^{jl}\omega_{ij}\omega_{kl}~$
$$\dt\vert{\omega}\vert^2=-2g^{in}g^{km}g^{jl}\omega_{ij}\omega_{kl}\dt g_{mn}
+2g^{ik}g^{jl}\omega_{ij}\dt\omega_{kl}$$
Using Lemma 1.4 and 1.7 we obtain
$$\eqalign{\dt\vert{\omega}\vert^2&=
4g^{in}g^{km}g^{jl}\omega_{ij}\omega_{kl}\eta^{st}H_sh_{tmn}
+2g^{ik}g^{jl}\omega_{ij}(\cd {H_l}k-\cd {H_k}l)\cr
&=4g^{km}\ud \omega nj\du \omega kj\eta^{st}H_sh_{tmn}
+2g^{ik}g^{jl}\omega_{ij}g^{pq}
\bigl(\dddd {\b R}{\underline p}qlk+
\cdd {\omega_{kl}}qp+\du {\omega}ps\dddd {\b R}sqlk\cr
&+\du {\omega}ls\dddd Rspkq+\du {\omega}ks\dddd Rspql
+\eta^{mn}\du {\omega}ns(h_{mkq}h_{slp}-h_{mlq}h_{skp})\bigr)\cr
&=4\omega^{nl}\ud {\omega}ml\eta^{st}H_s\ddd htmn+2\omega^{kl}
\dudd {\b R}{\underline p}plk+\Delta\vert\omega\vert^2-2\vert\cd
{\omega_{kl}}i\vert^2\cr
&+2 \omega^{kl}\du {\omega}ps \dudd {\b R}splk
+4\omega^{kl}\du{\omega}ls\dddu Rspkp+\omega^{kl}\eta^{mn}\du {\omega}ns(
\ddu hmkp\ddd hslp-\ddu hmlp\ddd hskp)\cr
}$$
Now it can be easily checked that each absolute term which depends
quadratically on $~\omega~$ can be bounded above by a term
$~c\vert\omega\vert^2~$ since $~L\times[0,T]~$ is compact, e.g. in
normal coordinates we see
that this is true for a term of the
form $~2\omega^{sl}\ud {\omega}mla_{sm}~$:
$$\eqalign{2\omega^{sl}\ud {\omega}mla_{sm}&=2\sum_{s,l}\bigl(\omega_{sl}
\sum_{m}\omega_{ml}a_{sm}\bigr)\le \vert\omega\vert^2+
\sum_{s,l}(\sum_{m}\omega_{ml}a_{sm})^2\cr
&\le \vert\omega\vert^2+n\sum_{s,l}\sum_m(\omega_{ml})^2(a_{sm})^2
\le (1+n\vert a\vert^2)\vert\omega\vert^2\cr}$$
These constants
depend on the dimension $~n~$, on $~\eta_{ij}~$,$~\ddd hijk~$, $~T~$
and the Riemann curvatures.
This gives
$$\dt\vert\omega\vert^2\le\Delta\vert\omega\vert^2+c\vert\omega\vert^2+
2\omega^{kl}\dudd {\b R}{\underline p}plk$$
The only thing which we have to prove is that the last term in the
above inequality can be bounded by a term of the form $~c\vert\omega\vert^2~$.
This will follow from the Ricci flatness of the ambient space. Let us assume
that $~v_i~$ is an orthonormal basis of eigenvectors for $~\eta_{ij}~$
and that, at a fixed point $~p\in L~$ we have chosen normal coordinates
for $~L~$ such that $~e_i=v_i~$ and that we have also chosen normal
coordinates for $~M~$ such that the corresponding coordinate frame is given
by $~e_i,~N(e_i)~$. Then the metric
$~\b g^{\alpha\beta}~$ at this point is given by
$$\b g^{\alpha\beta}=\pmatrix g^{ij}& 0\\ 0&\eta^{ij} \endpmatrix
=\pmatrix 1       &      &      &              &       & \\
                  &\ddots&      &              &       & \\
                  &      &1     &              &       & \\
                  &      &      &{1\over 1-a_1}&       & \\
                  &      &      &              &\ddots & \\
                  &      &      &              &       &{1\over 1-a_n}
\endpmatrix $$
where $~a_i:=\sum_{l=1}^n \omega _{il}^2~$.
For the Ricci curvature we obtain
$$0=\b Ric(Y,X)=\sum_{i=1}^n\langle\b R(X,e_i)e_i,Y\rangle+\sum_{i=1}^n
{1\over 1-a_i}\langle\b R(X,N(e_i))N(e_i),Y\rangle$$
Using the K\"ahler identity we see that this is equal to
$$\eqalign{0&=\sum_{i=1}^n\langle\b R(X,e_i)J(e_i),J(Y)\rangle+\sum_{i=1}^n
{1\over 1-a_i}\langle\b R(X,N(e_i))J(N(e_i)),J(Y)\rangle\cr
&=\sum_{i=1}^n\langle\b R(X,e_i)N(e_i),J(Y)\rangle
+\sum_{i,l=1}^n\omega_{il}\langle\b R(X,e_i)e_l,J(Y)\rangle\cr
&-\sum_{i=1}^n\langle\b R(X,N(e_i))e_i,J(Y)\rangle
-\sum_{i,l=1}^n\omega_{il}\langle\b R(X,N(e_i))J(e_l),J(Y)\rangle\cr
&+\sum_{i=1}^n{a_i\over 1-a_i}\langle\b R(X,N(e_i))J(N(e_i)),J(Y)\rangle\cr}$$
Now we use the Bianchi identity and see that
$$\eqalign{\sum_{i=1}^n\langle\b R(e_i,N(e_i))X,J(Y)\rangle&=
\sum_{i,m=1}^n\omega_{im}\langle\b R(X,e_i)e_m,J(Y)\rangle\cr
&-\sum_{i,m=1}^n\omega_{im}\langle\b R(X,N(e_i))J(e_m),J(Y)\rangle\cr
&+\sum_{i=1}^n{a_i\over 1-a_i}\langle\b R(X,N(e_i))J(N(e_i)),J(Y)\rangle\cr}$$
If we choose $~X=e_l,~Y=-J(e_k)~$ then we get
$$\eqalign{2\omega^{kl}\dudd {\b R}{\underline p}plk&=
2\sum_{i,k,l,m=1}^n\omega_{kl}\omega_{im}\dddd {\b R}kmli
-2\sum_{i,k,l,m=1}^n\omega_{kl}
\omega_{im}\langle\b R(e_l,N(e_i))J(e_m),e_k\rangle\cr
&+2\sum_{i,k,l=1}^n\omega_{kl}{a_i\over 1-a_i}\langle\b
R(e_l,N(e_i))J(N(e_i)),e_k\rangle\cr}$$
and since $~a_i\le\vert\omega\vert^2~$ and $~{1\over 1-a_i}~$ is bounded as
long as $~N~$ is an isomorphism, we can bound each of these terms by a
constant multiple of $~\vert\omega\vert^2~$. That
eventually proves Proposition 1.8
\enddemo

\proclaim{Theorem 1.9}
If $~L~$ is a closed manifold,
smoothly immersed as a Lagrangian submanifold
into a Ricci-flat Calabi-Yau manifold, then equation $(*)$ has a unique,
smooth solution at least on some short time interval and as long as the
solution
exists all submanifolds $~L_t~$ are Lagrangian submanifolds.
\endproclaim

\demo{Proof}
Existence and uniqueness follows from Proposition 1.6. To prove that this flow
preserves the Lagrangian structure we fix a time interval $~[0,T]~$ and look at
the function $~f_{\epsilon}:=\vert\omega\vert^2-\epsilon e^{2ct}~$, where
$~\epsilon~$ is an arbitrary, positive constant and
$~c~$ is as in Proposition 1.8.
Obviously $~f_{\epsilon}(0)=-\epsilon<0~$. On the other hand we obtain from
1.8 that
$$\dt f_{\epsilon}\le\Delta f_{\epsilon}+cf_{\epsilon}-c\epsilon e^{2ct}<
\Delta f_{\epsilon}+cf_{\epsilon}$$
Thus by the parabolic maximum principle we get that $~f_{\epsilon}<0~$ for
all $~t\in[0,T]~$ and all positive $~\epsilon~$. Then we let $~\epsilon~$
tend to zero and obtain the inequality
$\vert\omega\vert^2\le 0$ on $~[0,T]~$. Since also
$~\vert\omega\vert^2\ge 0~$ this shows that $~\vert\omega\vert^2~$
has to vanish identically and thus
$~\omega\equiv 0~$. This proves the theorem.
\enddemo

\head                 2. The Lagrangian case
\endhead

We turn our attention to the Lagrangian submanifolds. Since
$~\omega\equiv0,~\eta_{ij}=g_{ij},~N=J~$
we can restate the equations in 1.1 and 1.3 which become very beautiful
\proclaim {Proposition 2.1}
\formula{(a)}{h_{ijk}=h_{jik}=h_{jki}}
\formula{(b)}{\pddu F{\gamma}xkxj-\con nkj\pdu F{\gamma}xn+
\conbar {\gamma}{\alpha}{\beta}\pdu F{\alpha}xk\pdu F{\beta}xj=
-g^{mn}h_{mkj}\ud J{\gamma}{\beta}\pdu F{\beta}xn}
\formula{(c)}{\cd {\ddd hkji}l-\cd {\ddd hlji}k=\dddd {\b R}lkj{\underline i}}
\formula {(d)}{\dddd Rijkl=\dddd {\b R}ijkl+g^{mn}(h_{mik}h_{njl}
-h_{mil}h_{njk})}
\endproclaim

These equations are true for any Lagrangian submanifold. If in addition $~L~$
is immersed into a Ricci-flat Calabi-Yau manifold then we also have the
important identity
\formula {(e)}{dH=0}

In the forthcoming we will set $~\nu_s:=\nu_s^\alpha{\partial\over
\partial y^\alpha}:=\ud J{\alpha}{\beta}e_s^\beta
{\partial\over\partial y^\alpha}~$ and $~e_{ij}^\alpha:=
{\partial^2 F^\alpha\over\partial x^i\partial x^j}~$.

\proclaim{Proposition 2.2}
\formula{(a)}{\ddd hjkl=\dd {\b\omega}\alpha\beta e_{jk}^\alpha e_l^\beta
+\dd {\b\omega}\delta\gamma\conbar\delta\alpha\beta e_j^\alpha e_k^\beta
e_l^\gamma}
\formula{(b)}{{\partial \nu_s^\gamma\over\partial x^i}=\con lis\nu_l^\gamma
+\ddu hsil e_l^\gamma-e_i^\alpha\nu_s^\beta\conbar \gamma\alpha\beta}
\endproclaim

\demo {Proof}
(a) can be obtained by multiplying 2.1 (b) with $~\b g_{\gamma\delta}
\nu_l^\delta~$. For (b) we calculate
$$\eqalign{
\cdb {\nu_s}{e_i}&=\langle \cdb {\nu_s}{e_i},\nu_k\rangle\nu_l g^{kl}
+\langle \cdb {\nu_s}{e_i}, e_k\rangle e_l g^{kl}\cr
&=\langle \cdb {e_s}{e_i},e_k\rangle \nu_lg^{kl}-\langle \nu_s,
\cdb {e_k}{e_i}\rangle e_l g^{kl}\cr
&=\langle \cd {e_s}{e_i},e_k\rangle \nu_lg^{kl}+\ddd hsik e_l g^{kl}\cr
&=\con lis \nu_l+\ddu hsil e_l\cr}$$
On the other hand we have
$$\cdb {\nu_s}{e_i}=e_i^\alpha
\cdb {\bigl(\nu_s^\beta{\partial\over \partial y^\beta}\bigr)}{\alpha}=
\bigl({\partial \nu_s^\gamma\over \partial x^i}+e_i^\alpha\nu_s^\beta
\conbar \gamma\alpha\beta\bigr){\partial\over\partial y^\gamma}$$
Both equations together imply (b).
\enddemo

The evolution equation now takes the form
\formula {(*)}{\dt F^{\alpha}=-g^{mn}H_m\ud J{\alpha}{\beta}\pdu F{\beta}xn
=-H^n\nu_n^\alpha}
We also want to restate Lemma 1.7

\proclaim{Lemma 2.3}
\formula{(a)}{\dt g_{ij}=-2H^lh_{lij}}
\formula {(b)}{\dt d\mu=-H^lH_ld\mu}
\formula {(c)}{\dt \omega_{ij}=0}
\endproclaim

Now we derive the evolution equation for $~\ddd hjkl~$. The best way
to do this is to assume that we have chosen
double normal coordinates
at a given fixed point $~p\in L~$, i.e. normal coordinates for $~L~$ and $~M~$.
With Proposition 2.2 (a) and $(*)$ we obtain
$$\dt \ddd hjkl=-\b \omega_{\alpha\beta}\pddu {}{}xjxk(
H^n\nu_n^\alpha)e_l^\beta-\b\omega_{\alpha\beta}e_{jk}^\alpha
{\partial\over\partial x^l}(H^n\nu_n^\beta)
-\b\omega_{\delta\gamma}\conbar \delta\alpha{\beta,\epsilon}H^n\nu_n^\epsilon
e_j^\alpha e_k^\beta e_l^\gamma$$
since all other derivatives vanish in normal coordinates.
Let us denote these three terms by $ A,~B~$ and $~C~$. For $~A~$ we get
$$-\b \omega_{\alpha\beta}\pddu {}{}xjxk( H^n\nu_n^\alpha)e_l^\beta
=\uddd Hn,jk g_{ln}+\b\omega_{\beta\alpha}(\udd Hn,j\nu_{n,k}^\alpha
+\udd Hn,k\nu_{n,j}^\alpha)e_l^\beta+\b\omega_{\beta\alpha}H^n\nu_{n,jk}^
\alpha e_l^\beta$$
Using 2.2 (b) we see that this simplifies to
$$\eqalign{A&=
\uddd Hn,jk g_{ln}+\b\omega_{\beta\alpha}H^n\nu_{n,jk}^\alpha e_l^\beta\cr
&=\uddd Hn,jk g_{ln}+\b\omega_{\beta\alpha}H^ne_l^\beta(\con mjn
\nu_m^\alpha+\ddu hnjm e_m^\alpha-e_j^\gamma\nu_n^\delta
\conbar\alpha\gamma\delta)_{,k}\cr
&=\uddd Hn,jk g_{ln}+\b\omega_{\beta\alpha}H^ne_l^\beta(
\con mj{n,k}\nu_m^\alpha-\ddu hnjm\udd hsmk\nu_s^\alpha-\conbar \alpha\gamma
{\delta,\eta}e_k^\eta e_j^\gamma \nu_n^\delta)\cr}$$
again since all other derivatives vanish and because of $~\b\omega_{\alpha
\beta}e_i^\alpha e_j^\beta=0~$. Finally
$$\eqalign{A&=\uddd Hn,jk g_{ln}+g_{lm}\con mj{n,k}H^n-H^n\ddu hnjm\ddd hlmk-
\b\omega_{\beta\alpha}H^n\conbar \alpha\gamma{\delta,\eta}e_l^\beta
\nu_n^\delta e_k^\eta e_j^\gamma\cr
&=(\udd Hn,j+\con njmH^m)_{,k}g_{ln}
-H^n\ddu hnjm\ddd hlmk-
\b\omega_{\beta\alpha}H^n\conbar \alpha\gamma{\delta,\eta}e_l^\beta
\nu_n^\delta e_k^\eta e_j^\gamma\cr
&=(\nabla_jH^n)_{,k}g_{ln}
-H^n\ddu hnjm\ddd hlmk-
\b\omega_{\beta\alpha}H^n\conbar \alpha\gamma{\delta,\eta}e_l^\beta
\nu_n^\delta e_k^\eta e_j^\gamma\cr
&=\cdd {H_l}kj-H^n\ddu hnjm\ddd hlmk-
\b\omega_{\beta\alpha}H^n\conbar \alpha\gamma{\delta,\eta}e_l^\beta
\nu_n^\delta e_k^\eta e_j^\gamma\cr}$$
since in normal coordinates $~\cdd {H^n}kj=(\cd {H^n}j)_{,k}~$.
For $~B~$ we obtain
$$\eqalign{B&=-\b\omega_{\alpha\beta}e_{jk}^\alpha
{\partial\over\partial x^l}(H^n\nu_n^\beta)=\b\omega_{\alpha\beta}
\udd hsjk\nu_s^\alpha(\udd Hn,l\nu_n^\beta+H^n\ddu hnlm e_m^\beta)
=-H^n\udd hskj\ddd hnls\cr}$$
where we have used 2.1 (b) and 2.2 (b).
Now $~\b\omega_{\beta\alpha}e_l^\beta=\b g_{\alpha\beta}\nu_l^\beta~$.
Combining
$~A,~B~$ and $~C~$ gives
$$\eqalign{A+B+C&=
\cdd {H_l}kj-H^n(\ddu hnjm\ddd hmkl+\ddu hnlm\ddd hmkj)\cr
&-H^n(\b g_{\alpha\beta}\conbar\alpha\gamma{\delta,\eta}\nu_l^\beta
\nu_n^\delta e_k^\eta e_j^\gamma-\b g_{\delta\gamma}\conbar\delta\alpha
{\beta,\epsilon}\nu_n^\epsilon\nu_l^\gamma e_j^\alpha e_k^\beta)\cr}$$
After rearranging indices in the last term we see that this term equals
$$-H^n(\b g_{\alpha\beta}\conbar\alpha\gamma{\delta,\eta}\nu_l^\beta
\nu_n^\delta e_k^\eta e_j^\gamma-\b g_{\delta\gamma}\conbar\delta\alpha
{\beta,\epsilon}\nu_n^\epsilon\nu_l^\gamma e_j^\alpha e_k^\beta)=
H^n\b g_{\alpha\beta}(\conbar\alpha\gamma{\eta,\delta}-\conbar
\alpha\gamma{\delta,\eta})\nu_n^\delta\nu_l^\beta e_k^\eta e_j^\gamma$$
But since we have chosen normal coordinates this is exactly
$$H^n\b g_{\alpha\beta}(\conbar\alpha\gamma{\eta,\delta}-\conbar
\alpha\gamma{\delta,\eta})\nu_n^\delta\nu_l^\beta e_k^\eta e_j^\gamma
=H^n\dddd {\b R}\eta\delta\gamma\beta
\nu_n^\delta\nu_l^\beta e_k^\eta e_j^\gamma=H^n\dddd {\b R}{\underline n}
k{\underline l}j$$
Therefore we obtain the result

\proclaim {Lemma 2.4}
\formula{}{\dt \ddd hjkl=\cdd {H_l}kj-H^n(\ddu hnjm\ddd hmkl
+\ddu hnlm\ddd hmkj)+H^n\dddd {\b R}{\underline n}
k{\underline l}j}
\endproclaim

Note that the RHS of 2.4 is a symmetric three tensor since this is true for
$~h_{jkl}~$. This can also be proved directly by using the rule for
interchanging derivatives, the Gauss curvature equations, the first Bianchi
identity and the K\"ahler identity.

If we assume that $~L~$ is being
deformed by a different 1-form $~\theta~$ then it
can be easily checked that the
same calculations as above give the general evolution
equations

\proclaim {Lemma 2.5}
\formula{(a)}{\dt g_{ij}=-2\theta^lh_{lij}}
\formula {(b)}{\dt d\mu=-\theta^lH_ld\mu}
\formula {(c)}{\dt \omega=d\theta}
\formula {(d)}{\dt \ddd hjkl=\cdd {\theta_l}kj-\theta^n(\ddu hnjm\ddd hmkl
+\ddu hnlm\ddd hmkj)+\theta^n\dddd {\b R}{\underline n}
k{\underline l}j}
\endproclaim

Let $~d^\dagger~$ denote the negative adjoint to $~d~$, i.e. $~d^\dagger\theta
=\nabla_i\theta^i=g^{ij}\nabla_i\theta_j~$.
We want to calculate the time derivative of the mean curvature form

\proclaim{Lemma 2.6}
\formula {}{\dt H=dd^\dagger\theta}
\endproclaim

\demo {Proof}
Using 2.5 (a) and (d) we see that
$$\dt H_k=\dt(g^{jl}h_{jkl})=2g^{jm}g^{ln}\theta^s h_{smn}h_{jkl}+
\nabla_k d^\dagger \theta-2\theta^n\duu hnlm\ddd hmkl+\theta^n\dddu {\b R}
{\underline n}k{\underline l}l$$
The last term is a Ricci curvature and vanishes and two terms cancel
which gives
$$\dt H_k=\nabla _kd^\dagger\theta~.$$
\enddemo

$L~$ is called special Lagrangian, if it is Lagrangian and calibrated with
respect to the real part of the
Calabi-Yau form. 
Since calibrated manifolds are minimal
we obtain as a direct consequence of 2.5 (c) and 2.6:

\proclaim{Corollary 2.7}
If $~L~$ is a family of special Lagrangian surfaces evolving under $~\theta~$
then $~\theta~$ has to be closed and coclosed, i.e. harmonic.
\endproclaim

Let us finally prove the following identity which might be useful for
further investigations of the mean curvature flow.

\proclaim {Proposition 2.8}
$$\eqalign{
\cdd {\ddd hrsk}ij&
=\cdd {\ddd hijk}rs
+\cdb {\dddd {\b R}jrs{\underline k}}i
+\cdb {\dddd {\b R}isj{\underline k}}r\cr&
+\udd hnsk\dddd {\b R}njri
+\udd hnjk\dddd {\b R}nsri
+\udd hnjs\dddd {\b R}nkri\cr&
+\udd hnis\dddd {\b R}nkrj
+\udd hnik\dddd {\b R}nsrj
+\udd hnrj\dddd {\b R}nksi
+\udd hnrk\dddd {\b R}njsi\cr&
-\udd hnir\dddd {\b R}{\underline n}j{\underline k}s
+\udd hnij\dddd {\b R}{\underline n}r{\underline k}s
-\udd hnrs\dddd {\b R}{\underline n}i{\underline k}j
+\udd hnri\dddd {\b R}{\underline n}s{\underline k}j\cr&
+\udd hnsk(\ddu hnrm\ddd hmji-\ddu hnim\ddd hmjr)
+\udd hnjk(\ddu hnrm\ddd hmsi-\ddu hnim\ddd hmsr)\cr&
+\udd hnjs(\ddu hnrm\ddd hmki-\ddu hnim\ddd hmkr)\cr}$$
\endproclaim

\demo {Proof}
Using 2.1 (c) and the rule for interchanging derivatives we obtain
$$\eqalign{
\cdd {\ddd hrsk}ij&
=\cdd {\ddd hjsk}ri
+\cd {\dddd {\b R}jrs{\underline k}}i
+\udd hnsk\dddd { R}njri
+\udd hnjk\dddd { R}nsri
+\udd hnjs\dddd { R}nkri\cr&
=\cdd {\ddd hjik}rs
+\cd {\dddd {\b R}jrs{\underline k}}i
+\cd {\dddd {\b R}isj{\underline k}}r\cr&
+\udd hnsk\dddd { R}njri
+\udd hnjk\dddd { R}nsri
+\udd hnjs\dddd { R}nkri\cr}$$
Then the result follows from the Gauss curvature equations and the fact that
$$\eqalign{
\cdb {\dddd {\b R}jrs{\underline k}}i&=
\cd {\dddd {\b R}jrs{\underline k}}i
-\udd hnij\dddd {\b R}{\underline n}r{\underline k}s
+\udd hnri\dddd {\b R}{\underline n}j{\underline k}s\cr&
+\udd hnis\dddd {\b R}nkjr
-\udd hnik\dddd {\b R}nsrj\cr}$$
\enddemo

\proclaim {Theorem 2.9}
Assume that $~L_t=F_t(L)~$ is compact, orientable and
evolves under the mean curvature flow and that
$~x_0~$ is an arbitrary but fixed point on $~L~$.
\item{(a)} There exists a unique smooth family of functions $~\phi_t~$,
smoothly depending on time such that
$$H_t(x)=H_0(x)+d\bigl(\int_0^t\Delta_{\tau}\phi_{\tau}(x)d\tau\bigr)$$
$$\Delta_t\bigl(\phi_t-\int_0^t\Delta_{\tau}\phi_{\tau}(x)d\tau\bigr)=
d^\dagger_t H_0$$
$$\phi_t(x_0)=0$$
in particular the cohomology class of $~H~$ does not change.
\smallskip\noindent
\item{(b)} If the first Betti number of $~L~$ vanishes, then there exists a
unique smooth family of functions $~\phi_t~$ such that
$$H_t=d\phi_t$$
$$\dt\phi_t=\Delta_t\phi_t$$
$$\min_{L}\phi_0=0$$
\endproclaim

\demo {Proof}
Define the form $~\t
H_t(x):=H_0(x)+d\bigl(\int_0^td^\dagger_{\tau}H_{\tau}(x)d\tau\bigr)~$ where
we integrate pointwise. This form surely exists,
since $~d^\dagger_tH_t~$ is smooth. For the time derivative we obtain
$$\dt \t H_t=dd^\dagger_tH_t=\dt H_t$$ and since
$~\t H_0=H_0~$ we conclude
$~\t H_t=H_t~$. Now we use the decomposition theorem
and can express $~H_t~$ as a unique sum $~H_t=\psi_t+d\phi_t~$, where
$~d^\dagger_t\psi_t=0~$, $\phi_t(x_0)=0~$ and $~\psi_t,~\phi_t~$ are smooth.
Then $~d^\dagger_tH_t=d^\dagger_t d\phi_t=\Delta_t\phi_t~$. This proves (a).
(b) is a direct consequence of (a) since then the harmonic part $~\psi_t
\equiv 0~$ and therefore
$$H_t-H_0=d(\phi_t-\phi_0)
=d\bigl(\int_0^t\Delta_{\tau}\phi_{\tau}(x)d\tau\bigr)$$
This implies
$$d\bigl(\int_0^t\dt\phi_{\tau}-\Delta_{\tau}\phi_{\tau}(x)d\tau\bigr)=0$$
which means that there exists a smooth function $~f(t)~$ such that
$~\dt\phi_t-\Delta_t\phi_t=f(t)~$. Now define $~\t\phi_t:=\phi_t-\int_0^t
f(\tau)d\tau-min_{L}\phi_0~$. This function has all the desired properties.
If $~\phi_t,~\t\phi_t~$ are two functions with the same properties then
$~d(\phi_t-\t\phi_t)=0~$ and consequently there
exists a function $~f(t)~$ such that $~\phi_t=\t\phi_t+f(t)~$. Since
$~\dt(\phi_t-\t\phi_t)=\Delta_t(\phi_t-\t\phi_t)=0~$ we conclude that this
function is a constant $~c~$ which has to be zero since $~\min_{L}
\phi_0=\min_{L}\t\phi_0=0~$. This proves uniqueness.
\enddemo

\bigskip

\remark {Remarks and questions}
There are many interesting questions arising from this context. First of all
the above calculations show that the mean curvature form
of a Lagrangian submanifold is closed if the ambient space
is a Ricci-flat Calabi-Yau manifold. Therefore it makes sense to define the
mean curvature class $~[H]~$ as the cohomology class of the mean curvature
form. If $~L~$ is a minimal Lagrangian submanifold then clearly $~[H]=0~$.
The question is, if any Lagrangian submanifold with trivial
mean curvature class can be deformed into a minimal Lagrangian submanifold
of the same topological type. If $~L~$ is compact, orientable
and if the first Betti number of $~L~$ vanishes then $~H~$ is exact.
Are minimal Lagrangian submanifolds automatically special
Lagrangian? Under what conditions can one deform a Lagrangian submanifold
into a special Lagrangian? Is this possible if and only if the mean
curvature class of the initial surface vanishes?
We also note that the mean curvature vector has to vanish at least
in one point if the Lagrangian surface cannot admit an everywhere
nonzero, smooth tangent vectorfield because the complex structure
would map this vector field to a smooth and nowhere
vanishing tangent vector field. If the mean curvature form is exact and $~L~$
is closed then the
mean curvature has to vanish at least in two points, since then $~H=d\phi~$ and
the smooth function $~\phi~$ has at least one maximum and one minimum on $~L~$.

Another remarkable fact is that the pair (Lagrangian,
special Lagrangian) seems to be analogous to the pair (U(n) holonomy, SU(n)
holonomy). Since it has turned out in the past that the mean curvature
flow and the Ricci flow are always giving analogous results, one might ask
under what conditions does the Ricci flow or perhaps a different
flow preserve the holonomy of a metric. There are already some interesting
results with respect to that question \cite 1.

If one wants to deform a special Lagrangian manifold through special Lagrangian
manifolds then the corresponding 1-forms have to be closed and coclosed, i.e.
harmonic. This is sufficient and necessary and has been proven in \cite 6
(we only proved that this is necessary).
Unfortunately there is no canonical harmonic 1-form known which could be
used to deform such a surface. The most obvious approach would be
to define a coordinate
system on the moduli space of special Lagrangian manifolds by requiring that
the cohomology class does not change, thus giving something like geodesic
normal
coordinates (see \cite 8). However there are still great difficulties in
particular the existence of such a flow has not been proven.
Assuming that a special Lagrangian manifold moves under
a flow generated by a harmonic 1-form
one immediately gets a good insight in the structure of the evolution
equations by studying the corresponding evolution equations for the mean
curvature flow on Lagrangian submanifolds.
This is very analogous to the situation of a hypersurface evolving
in an ambient space. In that case the equations induced by the mean
curvature flow and those coming from different curvature functions are very
similar (see \cite 7)
and the equations for a Lagrangian surface flowing along
its mean curvature are also analogous to those coming from
a deformation of special Lagrangian
surfaces by harmonic 1-forms as we have already seen.
Thus the study of the mean curvature flow on
the moduli space of Lagrangian submanifolds in a Ricci-flat Calabi-Yau
manifold would help much in understanding what happens to special Lagrangian
surfaces.

Any smooth submanifold has a mean curvature vector and therefore it
would be possible to investigate the mean curvature flow in arbitrary
codimensions. But so far this has only been done
for hypersurfaces (e.g. \cite 5).
The reason is that a
hypersurface has a unique (up to the sign) unit normal vector. Given a basis
of the tangent space one can therefore identify this basis with a basis
of the normal space which in that case has dimension 1. For a Lagrangian
submanifold this can be achieved by using the complex structure $~J~$ and
this makes it
possible to compare the different second fundamental forms on $~L~$
which is a severe problem for a general submanifold.

There might be a greater number of canonical deformations. Given any smooth
family of functions $~f(t)~$
on $~L~$ one gets closed 1-forms $~\theta_t~$,
namely $~\theta_t=df_t~$. Obviously $~df_t~$ is closed independent of the
metric on $~L~$ and therefore independent of the immersion of $~L~$ into
$~M~$.
The corresponding
flow (if it exists) would preserve the Lagrangian structure since always
$~\dt \omega=d\theta_t~$.
On the other hand there is a great number of canonical,
smooth functions on $~L~$. These functions are the eigenfunctions of the
Laplace operator, i.e. $~\Delta f_i+\lambda_i f=0~$.
It would be interesting to
study the behaviour of $~L~$ under these flows, i.e. with $~\Delta_tf(t)
+\lambda_i(t)f(t)=0~$.
\endremark

\subhead Acknowledgements
\endsubhead
During the preparation of this paper the author has been a Post-Doc fellow
of the Alexander von Humboldt Foundation and would like
to thank S.T. Yau for his hospitality and his interest in this paper.
The author also thanks E. Zaslow for useful discussions.


\Refs
\ref\key 1
\by H.D. Cao
\paper Deformation of K\"ahler metrics to
K\"ahler-Einstein metrics on compact K\"ahler manifolds
\jour Invent. Math.
\vol 81\yr 1985\pages 359--372
\endref

\ref\key 2
\by R.S. Hamilton
\paper Three manifolds with positive Ricci curvature
\jour J. Differ. Geom.
\vol 23\yr 1986\pages  69--96
\endref

\ref\key 3
\by F.R. Harvey
\book Spinors and Calibrations
\publ Academic Press
\yr 1990
\endref

\ref\key 4
\by R. Harvey, H.B. Lawson
\paper Calibrated Geometries
\jour Acta Math.
\vol 148\yr 1982\pages 47--157
\endref

\ref\key 5
\by G. Huisken
\paper Contracting convex hypersurfaces in Riemannian
       manifolds by their mean curvature
\jour Invent. math.
\vol 84\yr 1986\pages 463--480
\endref

\ref\key 6
\by R. McLean
\paper Deformations of calibrated submanifolds
\jour Duke preprint
\yr 1996
\endref

\ref\key 7
\by K. Smoczyk
\paper Symmetric hypersurfaces in Riemannian manifolds
       contracting to Lie-groups by their mean curvature
\jour Calc. Var.\vol 4\yr 1996\pages 155--170
\endref

\ref\key 8
\by A. Strominger, S.T. Yau, E. Zaslow
\jour in preparation
\endref

\endRefs
\enddocument

\end